\documentclass{amsart}

\usepackage{graphicx}

\usepackage[all,cmtip]{xy}
\input{diagxy}

\usepackage{tikz}

\newtheorem{theorem}{Theorem}[section]
\newtheorem{lemma}[theorem]{Lemma}
\newtheorem{corollary}[theorem]{Corollary}

\theoremstyle{definition}
\newtheorem{definition}[theorem]{Definition}
\newtheorem{example}[theorem]{Example}

\theoremstyle{remark}
\newtheorem{remark}[theorem]{Remark}

\numberwithin{equation}{section}

\begin{document}
\title[Rational Homotopy Types with the Same Homotopy Lie Algebra]{A Moduli Space for Rational Homotopy Types with the Same Homotopy Lie Algebra}

\author{Matthew Zawodniak}

\date{12/7/15}
\thanks{The author gratefully acknowledges support from the RTG in Algebra, Algebraic Geometry, and Number Theory, at the University of Georgia (NSF grant DMS-1344994).} 

\maketitle

\section{Introduction}

Since Quillen proved his famous equivalences of homotopy categories in 1969, much work has been done towards classifying the rational homotopy types of simply-connected topological places. The majority of this work has focused on rational homotopy types with the same cohomology algebra. The models in this case were differential graded algebras and acted similarly to differential forms. These models were then used together with some deformation theory to describe a moduli space for all rational homotopy types with a given cohomology algebra. Indeed, this theory has been very well developed. However, there is another case to consider. That is, the collection of rational homotopy types with the same homotopy Lie algebra (same homotopy groups and Whitehead product structure). This case, arguably, is closer to the heart of homotopy theory, as it fixes the homotopy groups themselves and how they interact with each other. However, the Lie case has received less attention and is less developed than its cohomology counterpart.

The main purpose of this paper is to completely develop the theory for rational homotopy types of simply-connected topological spaces with the same homotopy Lie algebra. It will include some foundations of the theory as well as some new work. Often, previously-known results will be streamlined, reworded, or reproven to make them directly relevant to the results of this paper.

The first section of this paper will be a brief history of rational homotopy theory. The second will display the basics of the Lie algebra side of the theory. Both relevant constructions and relevant proofs will be provided here. In the final section, the deformation theory will be developed and the moduli space for rational homotopy types with a fixed homotopy Lie algebra will be defined and justified.

\section{A Brief History of Rational Homotopy Theory}

One of the main goals of homotopy theory has been to classify all homotopy types of simply connected topological spaces. Originally this meant integral homotopy types. While much work has been done on integral homotopy types \cite{blanc}, the integral homotopy groups of spheres, the basic building blocks of homotopy, have yet to be completely determined.

However, in 1951, Serre \cite{serre} showed that things worked out very nicely when looking at the \emph{rational} homotopy groups of spheres. Not only were the homotopy groups of spheres completed determined rationally, but they were also finite dimensional. As algebras, the homotopy groups of each sphere actually had only a single generator. This signaled that \emph{rationally}, spheres would behave quite well as the building blocks for a theory.

In 1969 Quillen made the field of rational homotopy theory even more desirable by showing that the category of rational homotopy types of simply connected topological spaces was equivalent to many other homotopy categories, including those of path-connected differential graded Lie algebras and simply-connected differential graded coalgebras. This led to the development of differential graded algebra (DGA) models \cite{sullivan}, which were later utilized in conjunction with some deformation theoretic ideas to make a moduli theory for simply-connected topological spaces with the same cohomology algebra. Shortly after this work on DGAs, many people began to work on the differential graded Lie algebra homotopy category. This paper will focus on differential graded Lie algebras as a means towards constructing a moduli space of rational homotopy types. Primary sources will be used whenever possible, though the text by F\'elix, Halperin, and Thomas gives a thorough development of the basics of rational homotopy theory.

\section{An Introduction to Lie Models}

\subsection{The Basics} We begin with a review of the basics.

\begin{definition}[Nijenhuis-Richardson \cite{NR}]{A \emph{graded Lie algebra}, or GLA, $L$ is a graded vector space together with a linear map $L\otimes L\to L$ given by $x\otimes y \mapsto [x,y]$ such that:

\begin{enumerate}
\item $[x,y]=(-1)^{|x||y|+1}[y,x]$, where $|x|$ is the degree of $x$.
\item $[x,[y,z]]=[[x,y],z]+(-1)^{|x||y|}[y,[x,z]]$, the graded Jacobi identity.

\

\end{enumerate}
A \emph{differential graded Lie algebra}, or DGLA, has the following additional properties:

\
\begin{itemize}
\item[(3)] a differential $d$ of degree $-1$.

\item[(4)] $d[x,y]=[dx,y]+(-1)^{|x|}[x,dy]$, i.e, $d$ is a graded derivation.

\

\end{itemize}}

\begin{remark} $(2)$ above basically says that $ad_x$ is a graded derivation, where $ad_x=[x,-]$.
\end{remark}

\begin{remark}This direction of the differential was Quillen's convention \cite{quillen} and is used throughout this paper.
\end{remark}
\end{definition}

The main motivation for looking at the category of differential graded Lie algebras comes from a result of Milnor-Moore, described below.

\begin{definition}
The \emph{homotopy Lie algebra} $P$ of a simply-connected topological space $X$ is the graded Lie algebra $P:=\pi_* (\Omega X)\otimes \mathbb{Q}$, where $\Omega X$ is the loop space on $X$.
\end{definition}

\begin{theorem}[Milnor-Moore {[M-M, p.263]}]
The homotopy Lie algebra $P$ as defined above is, in fact, a graded Lie algebra.
\end{theorem}

The bracket structure on $\pi_* (\Omega X)\otimes \mathbb{Q}$ is induced by the Whitehead product on $X$ and the long exact sequence of the loop space fibration on $X$, as seen in \cite{MM}.

\begin{remark} $\pi_* (\Omega X)\otimes \mathbb{Q}=\pi_{*+1} (X)\otimes \mathbb{Q}$, so the effect of using the loop space is basically just a shift of degree. This shift of degree is what allows the Whitehead product to respect the grading and induce a true bracket.
\end{remark}

Now that we see that the rational homotopy groups of a simply-connected topological space actually form a graded Lie algebra, it seems rather natural to look at the homotopy category of path-connected DGLAs to study rational homotopy types. In fact, in this category we are able to study rational homotopy types with the same homotopy Lie algebra.

\subsection{Models in the DGLA Homotopy Category}

\

In order to study this problem efficiently, there is a need for nice models or representations of the rational homotopy types of interest. Since we are working in the homotopy category of DGLAs, to us this means a kind of ``free" DGLA.

\begin{remark}
A graded free Lie algebra generated by the set $V$ will be denoted $\mathbb{L}_V$. A differential graded free Lie algebra is a graded free Lie algebra with a differential on it. Note that, in particular, this means that $\mathbb{L}_V$ is not necessarily free as a DGLA after being endowed with a differential.
\end{remark}

\begin{theorem}[Baues-Lemaire {[B-L, p.227-228]}, Niesendorfer {[N, p.448]}]
Given any DGLA $(L,d)$ there is a differential graded free Lie algebra, or DGFLA, $(\mathbb{L}_V,\bar{d})$ equipped with a quasi-isomorphism $(\mathbb{L}_V,\bar{d})\to (L,d)$. This is called the \emph{free Lie model} of $(L, d)$. Furthermore, there is a \emph{unique} free Lie model of $(L,d)$ in which $\bar{d}:V\to [\mathbb{L}_V,\mathbb{L}_V]$.

\end{theorem}

\begin{remark} If we take the homology of any of these models, we get the homotopy Lie algebra of the space we are modeling.
\end{remark}

\subsubsection{The Cellular Model} The existence and uniqueness result above for minimal models is definitely a requirement to allow the theory to move forward. However, it would be nice to have a more constructive approach available to us. In the case where $X$ is a finite-dimensional CW-complex, there is an incredibly straightforward construction for a corresponding Lie model.

\begin{theorem}[Neisendorfer {[N, p.458]}]
Given a simply-connected topological space $X$ with a CW-complex structure, we can build a free Lie model, called the \emph{cellular model}, as follows:
\begin{enumerate}
\item For each $n$-cell, $n\neq 0$, create a graded Lie algebra generator in degree $n-1$.
\item Define the differential $d$ by taking the attaching maps directly into this graded Lie algebra context. That is, if $e$ is a generator which represents a cell $U$ in the CW structure of $X$, then $d e$ is the representative of the image of $U$ under the attaching map.

\end{enumerate}
Furthermore, the homology of this DGLA yields the homotopy Lie algebra for $X$.
\end{theorem}

This theorem has many far-reaching results. One of the more obvious is that, given any free Lie model, we should be able to build a CW-complex of the same rational homotopy type that the free Lie model represents. However, before we get too far into the ramifications of the theorem, it would be useful to see an example of the cellular model in action.

\begin{example}[The cellular model of $S^2 $]

The most straightforward way to build $S^2$ as a CW-complex is with one $0$-cell and one $2$-cell with a trivial attaching map. Thus, we should have one generator in degree $1$ of the cellular model due to the shift of degree. The cellular model for $S^2$, then, is just
\begin{center}

\begin{tabular}{c c c c c c}
0 & $a$ & $[a,a]$ & $0$ & $0$ & ... \\
\hline
\hline
0 & 1 & 2 & 3 & 4 & ... \\
\end{tabular}
\end{center}
with the bottom row being the degrees. So, $S^2$ has $\pi_2$ and $\pi_3$ one-dimensional over $\mathbb{Q}$, with the rest $0$. This is consistent with the existence of the Hopf map $S^3\to S^2$.
\end{example}

The previous example, while simple, is perhaps too simple. The lack of nontrivial attaching maps means that only part of the cellular model is demonstrated. The next example is one step higher in complexity and will demonstrate the full potential of the cellular model.


\begin{example}[The cellular model of $\mathbb{C}P^2$]
$\mathbb{C}P^2$ has a standard CW decomposition of one $0$-cell, one $2$-cell, and one $4$-cell. The attaching map from the $4$-cell to the $2$-cell is actually given by the Hopf map. Thus, the cellular model looks like the following:


\begin{center}
\begin{tabular}{c c c c c c c}
0 & $a$ & $[a,a]$ & $b$ & $[b,a]$ & $[b,[a,a]]$ & ... \\
\hline
\hline
0 & 1 & 2 & 3 & 4 & 5 & ... \\
\end{tabular}
\end{center}
with the differential $d$ defined by $db=\frac{1}{2}[a,a]$. Note that, as in the $S^2$ model, $[a,a]$ corresponds to the image of the Hopf map. So we really are just encoding the data of the attaching map in the differential.
\end{example}

One of the properties of the cellular model, as noted above in the theorem, is that the homology of the cellular model is always the homotopy Lie algebra of the corresponding space. So, what is the homology of the DGLA above? $a$ is clearly a nonzero cycle in homology. However, $[a,a]$ is a boundary as $d(2b)=[a,a]$. The next nontrivial cycle is $[b,a]$, as $d[b,a]=\frac{1}{2}[[a,a],a]$ and $[[a,a],a]=0$ is a formal result of the graded Jacobi identity. Up through degree $5$, this is the only other nontrivial homology generator.

Now is a good time to note that these models are usually infinite dimensional! This means that there may be an infinite amount of homology generators above degree $5$. However, for our purposes knowledge of the homotopy Lie algebra up to degree $5$ will suffice.

\begin{center}
\begin{tabular}{c c c c c c c}
0 & $a$ & 0 & 0 & $[b,a]$ & 0 & ... \\
\hline
\hline
0 & 1 & 2 & 3 & 4 & 5 & ... \\
\end{tabular}
\end{center}
Above is the homotopy Lie algebra, up to degree $5$, of $\mathbb{C}P^2$. Note that $a$ and $[b,a]$ represent homology classes in this environment. However, to the nescient observer, it might look more like this:

\begin{center}
\begin{tabular}{c c c c c c c}
0 & $a$ & 0 & 0 & $x$ & 0 & ... \\
\hline
\hline
0 & 1 & 2 & 3 & 4 & 5 & ... \\
\end{tabular}
\end{center}
This is because without prior knowledge of the cellular model, ``$b$" has no meaning in the homotopy Lie algebra since it is not a cycle. Truly, the two graded Lie algebras are actually isomorphic. 

\subsubsection{The Coformal Model}

The homotopy Lie algebra $P$ serves a very important role in rational homotopy theory. In one sense, it is the simplest DGLA that has homology $P$, if given the zero differential. This makes $P$ serve as the crux of the corresponding rational homotopy types, as any other DGLAs with homology isomorphic to $P$ must be more complicated. The corresponding model and spaces of such a homotopy Lie algebra $P$ are of enough importance to warrant their own distinctive labels.

\begin{definition}[The Coformal Model]
The \emph{coformal model} is the unique minimal model quasi-isomorphic to its homotopy Lie algebra. A model is \emph{minimal} if its differential maps the underlying vector space $V$ to $[V,V]$. The space this model corresponds to is called the \emph{coformal space}. Note that sometimes $\pi$-formal is used in the literature. However, we use the term coformal because it is the analogue of formal models and spaces in the DGA case.
\end{definition}

\begin{remark}
This definition of coformal is similar but not identical to the definition in much of the literature (see \cite{NM} for one definition). However, it \emph{is} comparable and serves the same purpose. The main difference lies in the model considered. In most of the literature, the Quillen model is considered. However, any minimal model isomorphic to the Quillen model would also suffice. Even something like the cellular model of $S^2$, shown earlier, is a coformal model. For the rest of this section we will focus on a specific coformal model, the bigraded model.
\end{remark}

\subsubsection{The Bigraded Model}

The cellular model provided by Neisendorfer is a wonderfully straightforward tool to make models and find rational homotopy groups. However, it is a little too static to display all the rational homotopy types with its homotopy Lie algebra $P$. For this we need the \emph{bigraded model}. The bigraded model is just as algorithmic in construction, but big enough that we can adjust it to represent other, related rational homotopy types. While the process itself is rather straightforward, explaining it in generality is tedious. Rather than become overburdened with the technical details of the general case first, we proceed by way of example. The bigraded model for the homotopy Lie algebra of $\mathbb{C}P^2$ proceeds as follows:
\begin{example}[The bigraded model of the homotopy Lie algebra of $\mathbb{C}P^2$ \cite{oukili}]
\
\begin{enumerate}

\item Start with the homotopy Lie algebra.

\begin{center}
\begin{tabular}{c c c c c c c}
0 & $a$ & 0 & 0 & $x$ & 0 & ... \\
\hline
\hline
0 & 1 & 2 & 3 & 4 & 5 & ... \\
\end{tabular}
\end{center}
\item Choose a generating graded vector space $V$ for the homotopy Lie algebra.
\item Generate a free graded Lie algebra $\mathbb{L}_V$. This horizontal grading will be labeled ``topological degree".

\begin{center}
\begin{tabular}{c c c c c c c}
0 & $a$ & $[a,a]$ & 0 & $x$ & $[x,a]$ & ... \\
\hline
\hline
0 & 1 & 2 & 3 & 4 & 5 & ... \\
\end{tabular}
\end{center}

\item Look at the homology of $(\mathbb{L}_V,0)$. Clearly the homology of this DGFLA is not the same as the homotopy Lie algebra we started with. We need to adjust this DGFLA so that its homology is in fact the homotopy Lie algebra. To do this, we work from left to right and add generators where necessary. For instance, $[a,a]$ should be $0$ in homology. So, we add a generator $b$ in topological degree $3$ and resolution degree $1$ as shown below:

\begin{center}
\begin{tabular}{c || c c c c c c c}
1 &0 &  0  &   0     & $b$ &   &         &   \\
0 &0 & $a$ & $[a,a]$ & 0 & $x$ & $[x,a]$ & ... \\
\hline
\hline
&0 & 1 & 2 & 3 & 4 & 5 & ... \\
\end{tabular}
\end{center}
The vertical grading here is called the \emph{resolution degree}, so named as we are resolving the discrepancies between our model and the homotopy Lie algebra. For the rest of the paper, we will use upper indices for resolution degree and lower indices for topological degree. So, $\mathbb{L}^3_5$ will be the elements of topological degree $5$ and resolution degree $3$.

\

Now that we have our generator $b$, we define a differential $d$ by $db=[a,a]$, making this differential of bidegree $(-1,-1)$. Specifically, it reduces topological degree by $1$ and resolution degree by $1$. In this new model, $[a,a]$ is $0$ in homology as desired.

\item With this new generator added, generate the bigraded DGFLA.

\item Now we again scan resolution $0$ for any homology inconsistencies, starting from topological degree $0$. As $[x,a]$ is supposed to be $0$ in homology, we add in an generator $y$ in bidegree $(6,1)$ and define $dy=[x,a]$ to fix it.
\begin{center}
\begin{tabular}{c c c c c c c}
0 &  0  &   0     &  0  &    0     &       0     & $[b,b]$ \\
0 &  0  &   0     & $b$ & $[b,a]$  & $[b,[a,a]]$ & $y$  \\
0 & $a$ & $[a,a]$ & 0 & $x$ & $[x,a]$ & ... \\
\hline
\hline
0 & 1 & 2 & 3 & 4 & 5 & ... \\
\end{tabular}
\end{center}

\item Continue the consistency check until there are no unwanted non-boundary cycles in resolution degree $0$.

\item We would hope to be done, but adding in generators in higher resolution degree has added more unwanted elements with nonzero homology. For instance, $[b,a]$ is a cycle but not a boundary. So we add a generator $c$ in bidegree $(5,2)$ and define $dc=[b,a]$, then considering the resulting DGFLA.

\begin{center}
\begin{tabular}{c c c c c c c}
0 &  0  &   0     &  0  &    0     &       $c$     & $[b,b];[c,a]$ \\
0 &  0  &   0     & $b$ & $[b,a]$  & $[b,[a,a]]$ & $y$  \\
0 & $a$ & $[a,a]$ & 0 & $x$ & $[x,a]$ & ... \\
\hline
\hline
0 & 1 & 2 & 3 & 4 & 5 & ... \\
\end{tabular}
\end{center}
Note that $[b,[a,a]]$ is actually $0$ in homology from what we have already constructed, and thus we do not need to add a generator to kill it in homology.

\item Continue this process indefinitely or until it terminates. This will give you the bigraded model arising from the homotopy Lie algebra of $\mathbb{C}P^2$. This bigraded model will be quasi-isomorphic to $(P,0)$ as DGLAs.
\end{enumerate}
\end{example}

\

Now that we've seen an example it should be easier to follow the general construction of the bigraded model, outlined below:

\begin{enumerate}
\item Begin with a homotopy Lie algebra $P$ for some space.
\item Let the homotopy Lie algebra $P$ generate a free GLA $\mathbb{L}_V$, where $V$ is a graded sub-vector space of $P$ which maps isomorphically to $P/[P,P]$. This is a minimal generating set for $P$ as a graded Lie algebra.
\item Starting from topological degree $0$ and resolution degree $0$ and increasing topological degree, find the first element that does not show up in $P$ but is nonzero in homology.
\item Introduce a new generator one topological degree higher and one resolution degree higher than this element.
\item Define a differential $d$ mapping the new generator to the element which doesn't show up in $P$.
\item Generate a bigraded DGFLA with the old DGFLA and this added generator.
\item Repeat from step $3$ until all the elements in resolution degree $0$ have been adjusted as necessary.
\item Repeat from step $3$ on resolution degree $1$ now.
\item Keep successively adjusting each resolution degree until the homology of the bigraded DGFLA is the homotopy Lie algebra $P$ it was originally built from. Note that this process can be, and usually is, infinite.
\item Denote the end result of this process $(\mathcal{L}_P,d)$.

\end{enumerate}

\begin{remark}
Throughout the rest of this paper, we will use $(\mathcal{L}_P,d)$ to denote the bigraded model of $P$ or $(\mathcal{L},d)$ when $P$ is clear. The differential in the bigraded model reduces both resolution degree and topological degree by $1$, making this model bigraded as a DGLA.
\end{remark}

The bigraded model has many nice properties that make it a key structure in rational homotopy theory. The model is both free and minimal by construction.  In fact, the bigraded model is the unique minimal model mentioned in Theorem 3.8 when dealing with a coformal space. Of course, it is only unique up to isomorphisms of DGLAs, not bigraded DGLAs.

The second grading just allows us to keep track of how we are building the model. From the lens of the cellular model, the second grading is keeping track of sequences of cell attachments. Independent of viewpoint, the resolution degree plays no small role in characterizing rational homotopy types. The fact that the model is spread out in two dimensions makes it much easier to adjust in a controlled fashion.

Besides free, minimal, and resolution degree, the bigraded model has two other key properties.

\begin{lemma}
In the bigraded model, $\mathcal{L}_m^n=0$ when $ n-\frac{1}{2}m\geq 0$.
\end{lemma}

\emph{Proof:}
This follows directly from the construction of the bigraded model. Since the spaces we are dealing with are simply connected, $\mathcal{L}_0^0=0$. The first nonzero element can be in $\mathcal{L}_1^0$, call it $u$. If $u\neq 0$, then $u$ corresponds to a nonzero element of $P$, and thus is a nontrivial cycle. So, $\mathcal{L}_2^1=0$ and the first element in resolution degree $1$ must have at least topological degree $3$. Now, assume $\mathcal{L}_k^j=0$ for all $j\leq n$ and all $k\leq 2j$. We will show this holds true for resolution degree $n+1$. Now, as $\mathcal{L}^n_{\lceil \frac{1}{2}n \rceil}$ is the first nonzero dimension in resolution degree $n$, let $u\in \mathcal{L}^n_{\lceil \frac{1}{2}n \rceil}$. If $u$ is a generator of the bigraded DGLA, then it is not a cycle and therefore does not give rise to an element of resolution degree $n+1$. On the other hand, $u$ cannot be a bracket of other elements of $\mathcal{L}$ by the inductive hypothesis. So, the first boundary in resolution degree $n$ must be in topological degree $\lceil \frac{1}{2}n \rceil +1$, which implies that the first nonzero element in resolution degree $n+1$ must be of topological degree at least $\lceil \frac{1}{2}n \rceil +2$, as desired.
\qed

\begin{remark}
This diagonal of $0$s is visible in Ex. 3.15 above, and will be used to show local nilpotence in section $4$.
\end{remark}

\begin{lemma}
All nonzero homology classes of a bigraded model have representatives solely in resolution degree $0$.
\end{lemma}

\emph{Proof:}
This lemma follows directly from the construction of the bigraded model. The construction begins in resolution degree $0$ with the homotopy Lie algebra. Then we let this generate a free Lie algebra. Next we introduce elements in resolution degree $1$ that map by $d$ to elements of resolution degree $0$ that are not supposed to live in homology and add these elements as generators of our free lie algebra. Again, if we introduce any elements of $ker(d)$ we then add elements in the next resolution degree to bound them off. Thus, elements only exist in higher resolution degree if they bound off an element or are a boundary themselves.
\qed

Now we return to Ex. 3.12 and Ex. 3.11. Recall that we constructed both the cellular model for $\mathbb{C}P^2$

\begin{center}
\begin{tabular}{c c c c c c c}
0 & $a$ & $[a,a]$ & $b$ & $[b,a]$ & $[b,[a,a]]$ & ... \\
\hline
\hline
0 & 1 & 2 & 3 & 4 & 5 & ... \\
\end{tabular}
\end{center}
as well as the bigraded model of the homotopy Lie algebra $P$ of $\mathbb{C}P^2$, which we now know to be the bigraded model of the coformal space with homotopy Lie algebra $P$.

\begin{center}
\begin{tabular}{c c c c c c c}
0 &  0  &   0     &  0  &    0     &       $c$     & $[b,b];[c,a]$ \\
0 &  0  &   0     & $b$ & $[b,a]$  & $[b,[a,a]]$ & $y$  \\
0 & $a$ & $[a,a]$ & 0 & $x$ & $[x,a]$ & ... \\
\hline
\hline
0 & 1 & 2 & 3 & 4 & 5 & ... \\
\end{tabular}
\end{center}

A natural question, then, is if $\mathbb{C}P^2$ actually \emph{is} this coformal space, or at least quasi-isomorphic to it. Since both of these models are free and minimal, it would suffice to show that the two models are isomorphic as DGLAs. However, Neisendorfer and Miller showed that this was, in fact, not the case \cite{NM}.

\begin{example}[$\mathbb{C}P^2$ is not coformal]
While Neisendorfer and Miller prove this result in \cite{NM} by dualizing and using DGA results, this result can also be seen directly from the two models. As the models are more pertinent to the task at hand, we demonstrate the latter method below.

In this case, the fact that there is no isomorphism is rather straightforward. Suppose we attempted to build a DGLA isomorphism $\phi$ from the cellular model of $\mathbb{C}P^2$ to the coformal model. Then to start, we would be required to map $a$ to $a$. Then, to preserve the differential, the isomorphism $\phi$ would have to map $b$ to $b$. However, this leads to an inconsistency. To preserve the homology we would need $\phi [b,a]=x$, but to preserve the bracket we need $\phi [b,a]=[b,a]$. Thus such an isomorphism is impossible.

\end{example}

We saw in the previous example that the cellular model for $\mathbb{C}P^2$ and the coformal model are not isomorphic. By uniqueness of minimal models, then, these models represent two distinct rational homotopy types. However, they still have the same homotopy Lie algebra. So, one would expect the two to be related somehow. In fact, the observant reader may have already noticed that the cellular model injects into the bigraded model. However, the homology is not preserved. We could fix that if we were just allowed to changed the differential in the bigraded model by a little. For instance, $dc=[b,a]-x$ would make $x$ and $[b,a]$ represent the same class in homology and ``make" the two models equivalent. This slight change to the differential is an example of a \emph{perturbation}, which will be discussed thoroughly in the next section.

\section{The Deformation Theory}

Before describing perturbations, we first take a step back to examine {\it{why}} we study perturbations. Historically, this idea stems from Halperin-Stasheff and 

\noindent Schlessinger-Stasheff in the late 1970s \cite{H&S} \cite{SS2}. Halperin and Stasheff showed, in the fixed cohomology case, that every rational homotopy type with cohomology algebra $\mathcal{H}$ has a representative DGA model obtained by deforming a particular formal model for $\mathcal{H}$. These deformations are, in fact, suitable perturbations to the differential. One of the properties of said ``suitable" perturbations is that they are derivations. As derivations on a DGA form a DGLA, Schlessinger-Stasheff adopted the mindset that deformation problems in rational homotopy theory were governed by a ``controlling DGLA'' (as in \cite{SS1}). This also holds true in the fixed homotopy Lie algebra case, as we shall see in this section.

\subsection{The Perturbation Theorem}

Since the controlling DGLA will contain our perturbations, it is now the appropriate time to define what a perturbation is.

\begin{definition}
A \emph{perturbation} of a DGLA $(L,d)$ is a derivation $\tau$ of topological degree $-1$ such that $(d+\tau)^2=0$.
\end{definition}

This equation that perturbations must satisfy is often written $$D\tau + \frac{1}{2}[\tau,\tau]=0$$ and called the \emph{Maurer-Cartan} equation. The benefit of this form is that it is phrased in Lie algebra language. However, the differential $D$ is of the derivation Lie algebra, not of the model $(\mathcal{L},d)$ which is being perturbed.

\begin{definition}
The \emph{derivation Lie algebra} $Der \mathcal{L}$ is the graded Lie algebra of graded derivations of $\mathcal{L}$, where $Der_i \mathcal{L}=$ {derivations of degree $i$}. The bracket is given by the commutator, and the differential is given by $D:=ad_d(-)= [d,-]$, a simple bracketing with the differential on $\mathcal{L}$. As defined, $(Der\mathcal{L},ad_d)$ is a DGLA. 
\end{definition}

The scrutinizing reader can check for themselves that the Maurer-Cartan equation is in fact the same statement as $(d+\tau)^2=0$, i.e. that the perturbed differential $d+\tau$ is still a differential.

\

As mentioned at the beginning of this section, our interest in perturbations stems from the fact that, in the DGA case, one need only perturb the formal model to get all rational homotopy types with the same cohomology algebra. Thankfully, or perhaps as expected, the analogous statement holds in the DGLA case.

\begin{theorem}[Perturbation Theorem {[Ou, p.30]}]
Let $(L,\delta)$ be a differential graded Lie algebra with homotopy Lie algebra $P$. Let $(\mathcal{L},d)$ be the bigraded model of $P$. Then there exists a perturbation $\tau$ which reduces resolution degree by at least $2$ and a quasi-isomorphism $\pi:(\mathcal{L},d+\tau)\to (L,\delta)$. That is, $(\mathcal{L},d+\tau)$ is a DGLA model of $(L,\delta)$.
\end{theorem}

\emph{Proof:}(Adapted from \cite{oukili})

Let $\eta:P\to \mathcal{L}^0$ be a linear map for which $i\circ \eta = id$, where $i:(\mathcal{L},d)\to (P,0)$ induces an isomorphism on homology. We construct $\tau$ and $\pi$ by induction on resolution degree.

In resolution degree $0$, $\tau=0$ and we define $\pi$ so that the following diagram commutes:

$$\bfig
\morphism(0,0)|a|/{->}/<500,0>[\mathcal{L}^0`Z(L);\pi ]
\morphism(0,0)|b|/{->}/<250,-250>[\mathcal{L}^0`P;i]
\morphism(500,0)|b|/{->}/<-250,-250>[Z(L)`P;i']
\efig$$
where $i$ and $i'$ induce isomorphisms on homology and $Z(L)$ is the collection of cycles in $(L,\delta)$. Therefore, the class of $\pi(u)$, denoted $cl(\pi(u))$, is $i(u)$ for $u\in \mathcal{L}^0$ as a statement for $P$ naturally extends to $\mathcal{L}^0$.

On $\mathcal{L}^1$ we let $\tau =0$. Clearly $cl(\pi (du))=i(du)=0$ for $u$ a generator of $\mathcal{L}^1$. So, we can choose a linear map $\pi:\mathcal{L}^1\to L$ for which $\delta \pi (u)=\pi (du)$ by defining such a map on generators. Thus $\pi$ is a morphism of DGLAs.

We extend to $\mathcal{L}^2$ in the following way. Let $x\in \mathcal{L}^2$ be a generator of $\mathcal{L}$. We want to define $\tau x$ and $\pi x$ so that $\pi (d+\tau)x=\delta \pi x$. We will use the fact that $\pi dx$ was already defined in order to define $\tau$ and $\pi$ accordingly. First, as $dx$ is a cycle for $d+\tau$ ($\tau dx=0$ by a degree argument), $\pi dx$ must also be a cycle, as $\pi$ was defined to be a DGLA homomorphism through resolution degree $1$. This leaves us with two cases.

\underline{Case 1:} $\pi dx$ is a boundary for $\delta$.

Then we can define $\tau x:=0$ and $\pi x$ to be an element $a\in L$ such that $\delta a= \pi dx$.

\underline{Case 2:} $\pi dx$ is not a boundary for $\delta$.

Consider $\alpha=i'\pi dx$, the nonzero homology class of $\pi dx$, and define $\tau x = -\eta \alpha $. Then $\pi(d+\tau)x=\pi(dx- \eta \alpha)$, which is a boundary as $i'\pi dx=i'\pi \eta \alpha$ by construction. So, $\exists a\in L$ such that $\delta a = \pi(dx- \eta \alpha)$ and we define $\pi x:= a$.

Note that $(d+\tau)^2=d^2=0$ on $\mathcal{L}^{\leq 2}$, so $d+\tau$ is a differential as defined, and $\pi (d+\tau)=\delta \pi$ on $\mathcal{L}^{\leq 2}$, so $\pi$ as defined up to this stage is a DGLA homomorphism. Furthermore, $\eta \alpha \in \mathcal{L}^0$ and so $\tau$ reduced resolution degree by $2$. We now proceed to the inductive step.

\underline{Inductive Step:} Assume $d+\tau$ and $\pi$ have been defined through resolution degree $n$ for $n\geq 2$ such that $\tau$ reduces resolution degree by at least $2$. Let $u\in \mathcal{L}^{n+1}$ be a generator of $\mathcal{L}$. We consider, specifically, $\tau du$. Note that $\tau du=(d+\tau)du$, which means that $du$ is a boundary for $d+\tau$. So as $\pi$ is already defined through resolution degree $n$ and $\tau du$ is in resolution degree $\leq n-1$, $\tau du$ would have to be an element which was already a boundary for $d+\tau$. Name this element $w\in \mathcal{L}^{\leq n-1}$. Then $(d+\tau)w=\tau du$. So, if we define $\tau u=-w-\eta \alpha$, where $\alpha = i'\pi(du-w)$, then $(d+\tau)^2 u=0$ as desired. Furthermore, $i'\pi (d+\tau)u=0$, so $\exists a\in L$ such that $\delta a= \pi(d+\tau)u$ and we can define $\pi u:= a$.

This completes the construction of $\pi$ and $d+\tau$.
\qed

\

Haralambous proved a generalized version of this theorem in \cite{haralambous} using spectral sequences.

\begin{remark}
After perturbing the bigraded model it is no longer bigraded. Rather, it is only filtered graded as a differential graded Lie algebra, due to the fact that perturbations only respect the topological grading, as seen in the definition of the controlling DGLA below. 
\end{remark}

\subsection{The Controlling DGLA}

We know from the perturbation theorem that any rational homotopy type can be obtained by perturbing the bigraded model. However, as has been foreshadowed throughout this paper, not every perturbation corresponds to a rational homotopy type with the same homotopy Lie algebra. Furthermore, there can be instances where two perturbations correspond to the same rational homotopy type. Once these issues are taken care of we will be ready to create a moduli space for rational homotopy types with a fixed homotopy Lie algebra. (For a complete background on controlling DGLAs in deformation theory, see \cite{manetti})


\begin{example}[Two Isomorphic Bigraded Models]

Let $P$ be the following graded Lie algebra, generated by $a$, $b$, and $z$:

\begin{center}
\begin{tabular}{c c c c c c c}
0 & $a;z$ & $[a,a]$ & $b;[a,z]$ & $[b,a];[b,z]$ & $[b,[a,a]]$ & ... \\
\hline
\hline
0 & 1 & 2 & 3 & 4 & 5 & ... \\
\end{tabular}
\end{center}

The two bigraded models below correspond to this rational homotopy type. The first is built straight from $P$ and has homology equal to $P$, while the other would require an automorphism $a\to z$.

\begin{center}
\begin{tabular}{c c c c c}
0 & 0       & 0             &    $y$      & ... \\
0 & $a;z$ & $[a,a];[z,z]$ & $b;[a,z]$ & ... \\
\hline
\hline
0 & 1 & 2 & 3 & ... \\
\end{tabular}
\end{center}

{\center{with $dy=[z,z]$.}}

\

\begin{center}
\begin{tabular}{c c c c c}
0 & 0       & 0             &    $y$      & ... \\
0 & $a;z$ & $[a,a];[z,z]$ & $b;[a,z]$ & ... \\
\hline
\hline
0 & 1 & 2 & 3 & ... \\
\end{tabular}
\end{center}

{\center{with $dy=[a,a]$.}}

\

The fact that these two models are practically identical should not be surprising. After all, they represent the same rational homotopy type. With the added knowledge that the bigraded model can be built directly from $P$, we can choose the correct model, in this case the first one, as the model to represent the coformal space.

\end{example}

We can now define the collection of suitable perturbations which will completely characterize the rational homotopy types with a given homotopy Lie algebra.

\begin{definition} For a given bigraded model $(\mathcal{L}_P,d)$, we define $\Theta\subset Der\mathcal{L}_P$ to be the differential graded Lie subalgebra such that $\Theta_i=\{$derivations of degree $i$ that decrease resolution degree by more than $-i\}$. Specifically, any $\tau \in \Theta_{-1}$ must decrease resolution degree by more than $1$.

\end{definition}

\begin{remark} Note that the perturbation proposed at the end of section 3 is, in fact, an element of $\Theta_{-1}$.
\end{remark}

It has already been shown by Oukili, and later confirmed and generalized by Haralambous, that every rational homotopy type with a given homotopy Lie algebra can be obtained from perturbations of the bigraded coformal model of a simply connected topological space. Now, we must show that every perturbation $\tau \in \Theta_{-1}$ yields a DGLA with the same homotopy Lie algebra $P$.

\

Let $( \mathcal{L},d)$ be a DGLA with differential of degree $-1$ which is the minimal bigraded DGLA model for a coformal space. We define a splitting $\phi $ of $d$ as follows:

\

Consider $B_i \subset Z_i \subset L_i$, where $L_i$ are the chains of topological degree $i$ in $\mathcal{L}$, $Z_i$ the cycles, and $B_i$ the boundaries. For each $i$, we can decompose $L_i$ as a direct sum of vector spaces in the following way: First note that $Z_i=V_i\oplus B_i$ where $V_i$ is the complementary vector space to $B_i$ in $Z_i$. Similarly, $L_i=W_i\oplus Z_i$. So, we have $L_i=W_i\oplus V_i\oplus B_i \forall i$. The differential $d$ induces an isomorphism $d|_{W_{i+1}}:W_{i+1}\to B_i$ for each $i$. Let $\phi_B:B_i\to W_{i+1}$ be the inverse of this isomorphism. We extend $\phi_B$ to a $\phi_Z$ on $Z_i$ by setting $\phi_Z(V_i)=0$ and $\phi_Z|_{B_i}=\phi_B$ and extending linearly. We can further extend $\phi_Z$ to $\phi; L_i\to L_{i+1} $ by setting $\phi(W_i)=0$, $\phi|_{Z_i}=\phi_Z$, and extending linearly. This $\phi$ is a linear map which is a splitting of the differential $d$.

\

Now that we have a splitting, we can state the theorem.

\begin{theorem}
The map $f:(\mathcal{L},d)\to (\mathcal{L},d+\tau)$ defined by $f(x)=x+\tau \phi (x)$ is a bijection and induces an isomorphism of graded Lie algebras on homology, with $\tau$ and $\phi$ as defined above.
\end{theorem}

\emph{Proof:}

We first show bijectivity of $f$, and then show $f$ preserves homology.
\

\begin{itemize}
\item \underline{Injective:} 

It is enough to show that $f$ is injective on each $L_i$. So, let $x\in L_i$. We can view $x$ as a sum of elements in each resolution degree, $\sum_{j=0}^{k}x_j$ for some $k$, as $L_i$ is necessarily $0$ in resolution degrees $\geq i$. (Note, this bound is not sharp, though I could replace it with the sharp one at some point). Assume $f(x)=x+\tau \phi (x)=0$. Then as $\tau$ reduces resolution degree by at least $2$ and $\phi$ increases resolution degree by exactly $1$, $x_k=0$. Then $x=\sum_{j=0}^{k-1}x_j$. By the same argument, $x_{k-1}=0$. This continues all the way to $x_0=0$, and thus $x=0$.

\item \underline{Surjective:}

Let $x\in L_i$ and define $A:=\sum_{i=0}^\infty (-1)^i(\tau \phi)^i(x)$. Note that $A$ is actually a well-defined finite sum as $L_i$ consists of elements of resolution degree $\leq i$ and $\geq 0$, and $\tau \phi$ reduces resolution degree by at least $1$. Then $$f(A)=A+\tau \phi (A) = \sum_{i=0}^\infty (-1)^i(\tau \phi)^i(x)+\sum_{i=0}^\infty (-1)^i(\tau \phi)^{i+1}(x)$$
$$=\sum_{i=0}^\infty (-1)^i(\tau \phi)^i(x)+\sum_{i=1}^\infty (-1)^{i+1}(\tau \phi)^{i}(x)$$
$$=x+\sum_{i=1}^\infty (-1)^i(\tau \phi)^i(x)+\sum_{i=1}^\infty (-1)^{i+1}(\tau \phi)^{i}(x)=x$$

\

\item \underline{$f$ induces an isomorphism of graded Lie algebras:}

First we show $f:B_i^d\to B_i^{d+\tau}$. Let $y\in B_i^d$. Then $y=dz$ for $z=\phi y$, and $f(y)=y+\tau \phi(y)=dz+\tau z=(d+\tau)z\in B_i^{d+\tau}$. Thus $f$ induces a linear map $f_*H(\mathcal{L},d)\to H(\mathcal{L},d+\tau)$.

Now let $[x]\in H(\mathcal{L},d)$ and choose a representative $x\in Z_{i,0}^d$, which exists by Lemma 3.19. Then clearly $(d+\tau)x=0$ as $L_i$ is $0$ in negative resolution degree. Furthermore, $\tau \phi (x)=0$ for the same reason and so $f(x)=x\in Z_i^{d+\tau}$. 

\underline{Injective on Homology:} Assume $[x]\neq [0]$ in $H_i^d$. We want to show that $[x]\neq [0]$ in $H_i^{d+\tau}$. By way of contradiction, assume $x=(d+\tau)\sum_{i=1}^ny_i$ where $y_i$ has resolution degree $i$. We proceed by induction on $n$. If $n=1$, then $x=dy\in im(d)$, a contradiction. Now assume $x=(d+\tau)\sum_{i=1}^ky_i\implies x\in im(d)$ for $k<n$. Assume $x=(d+\tau)\sum_{i=1}^ny_i$. Then as $x\in Z_{i,0}$, $dy_n=0$. As $y_n$ has resolution degree $n$, by construction of the model $\exists z$ such that $dz=y_n$. Consider $\sum_{i=1}^ny_i-(d+\tau)z$. $(d+\tau)[\sum_{i=1}^ny_i-(d+\tau)z]=x$, but $\sum_{i=1}^ny_i-(d+\tau)z$ consists of elements of resolution degree $<n$. By our inductive hypothesis, $x\in im(d)$, a contradiction. Thus the linear map $f_*$ on homology is injective.

To show surjectivity on the homological level we need one more lemma.

\begin{lemma} If $0\neq[x]\in H_i(\mathcal{L},d+\tau)$, then $[x]$ has a representative in resolution degree $0$.

\end{lemma}

\emph{Proof:}
Let $\bar{x}$ be a representative for $[x]$ of least resolution degree. That is, $\bar{x}\in(\mathcal{L},d+\tau)_i^{\leq k}$ with $k$ the least possible. In other words, $\bar{x}=\sum_{j=0}^k\bar{x}^j$. In this case, we set $mrd(\bar{x})=k$. In other words, the maximum resolution degree piece of $\bar{x}$ lies in resolution degree $k$. By way of contradiction, assume that $k\geq 1$. Then as $(d+\tau)\bar{x}=0$, $d\bar{x}^k=0$. As the bigraded model has no nontrivial $d$-cycles in positive resolution degree by construction, $\bar{x}^k=dy$ for some $y$ by the structure of the minimal bigraded coformal Lie model. Now if $\tau y \neq 0$ then $\bar{x}\sim \sum_{j=0}^{k}\bar{x}^j-(d+\tau) y=\sum_{j=0}^{k-1}\bar{x}^j-\tau y$ in $H_i(\mathcal{L},d)$. But $mrd(\sum_{j=0}^{k-1}\bar{x}^j-\tau y)<k$ and we reach a contradiction of the minimality of $\bar{x}$. However, if $\tau y = 0$ then $(d+\tau)y=dy=x$ and so $[x]=0$, which is also a contradiction. So, $[x]$ must have a representative in resolution degree zero.
\qed

\

Now that we know that each element of homology has a representative in resolution degree $0$, we can proceed with surjectivity on the homological level.

\underline{Surjective on Homology:} Let $[x]\in H_i(\mathcal{L},d+\tau)$. By Lemma 4.9 above, we can choose a representative $\bar{x}$ of $[x]$ in resolution degree $0$. Furthermore, choose $\bar{x}$ such that $\bar{x}\in V_i^{d+\tau}$. As $\bar{x}\neq 0$, $\bar{x}\notin im(d+\tau)$. However, as $\bar{x}$ has resolution degree $0$, this is the same as saying $\bar{x}\notin im(d)$. So $[\bar{x}]\neq 0$ in $H_i(\mathcal{L}, d)$. As $f(\bar{x})=\bar{x}$, it follows that $f_*[\bar{x}]=[\bar{x}]=[x]$ as desired.

As $f$ was a bijection we see that the induced map is an isomorphism of graded vector spaces. As $f$ is the identity on $Z_{i,0}$ and thus induces the identity map on homology, the bracket structure is also preserved by $f_*$. It follows that $f_*$ induces an isomorphism of graded Lie algebras. \qed
\end{itemize}

We have shown that any perturbation $\tau\in \Theta_{-1}$ of the bigraded coformal model yields a filtered model whose homotopy Lie algebra is isomorphic to the homotopy Lie algebra of the coformal model, as desired.

One thing that we have not addressed is the dependence of $f$ on $\tau$. Technically, the isomorphisms $H(\mathcal{L},d+\tau)\to P$ for each $\tau$ and the automorphisms of $P$ itself can also cause rational homotopy types to be counted multiple times. For this reason, the literature often refers to a triple $(\mathcal{L},d+\tau,i_\tau)$ when considering the perturbation theory. However, one can apply an automorphism to $\mathcal{L}$ which causes all the isomorphisms to be identical in the following way. Let $V$ be a set of generators for $P$, and let $I:P\to (\mathcal{L},d+\tau)$ be a splitting of the homology map. Then $\phi \circ I (V)$ is a set of generators for $P$ as $\phi$ is a weak equivalence. So, the map $\phi \circ I (v) \mapsto I(v)$ induces an isomorphism of $P$, and thus an isomorphism of $(\mathcal{L}_P, d)$ by construction.

\begin{definition} A \emph{quasi-isomorphism of triples} $\phi:(\mathcal{L},d+\tau,i_\tau)\to (\mathcal{L},d+\tau',i_{\tau'})$ is a quasi-isomorphism $\phi:(\mathcal{L},d+\tau)\to (\mathcal{L},d+\tau')$ such that the following diagram commutes:

$$
\bfig
\morphism(0,0)|a|/{->}/<1000,0>[H(\mathcal{L},d+\tau)`H(\mathcal{L},d+\tau');\phi_* ]
\morphism(0,0)|a|/{->}/<0,-500>[H(\mathcal{L},d+\tau)`P;i_\tau]
\morphism(1000,0)|a|/{->}/<0,-500>[H(\mathcal{L},d+\tau')`P;i_{\tau'}]
\morphism(0,-500)|a|/{->}/<1000,0>[P`P;id]

\efig
$$

\end{definition}

\begin{remark} As a result of this theorem, $\Theta$ can be identified as the ``controlling DGLA" of the deformation theory here.
\end{remark}

\subsection{The Moduli Space Construction}

$\Theta$ as defined may include all of the suitable perturbations, but many of these perturbations can represent the same rational homotopy type. It would be preferable if such perturbations could be identified as ``equivalent" somehow. One way to do this is through what is called a \emph{gauge action}.
The gauge action comes from the standard action of $exp(\Theta_{0})$ on $\Theta_{-1}$. That is, $exp(ad_\theta)$ acts on $\Theta_{-1}$ for each $\theta \in \Theta_0$.

\begin{remark} The results of this paper depend rely heavily on not only the gauge action itself, but the interaction between the gauge action as defined and the bigraded model. Specifically, the gauge action is \emph{locally nilpotent} on the bigraded model (the curious reader can refer back to Example 3.15 to see that there is a line of slope $\frac{1}{2}$ on and above which lie only $0$). Local nilpotence implies two essential properties of the gauge action. First, the action is invertible, which will allow the statements below to generalize to arbitrary chains of weak equivalences regardless of direction. Second, the action preserves the bracket structure, which is essential to the following theorem.
\end{remark}

\begin{theorem}
There exists a quasi-isomorphism of triples $\phi :(\mathcal{L},d+\tau, i_\tau)\to (\mathcal{L}, d+\tau',i_{\tau'}) $ if and only if there exists a $\theta \in \Theta_0$ such that $ (d+\tau') = exp(ad_\theta )(d+\tau)$.

\end{theorem}

\emph{Proof:}
\underline{$\impliedby$:}

\begin{description}
\item[a)] Let $\theta\in\Theta_0\subset Der_0L$ and $\tau\in \Theta_{-1}$. Consider $exp(ad_\theta)(d+\tau)=d+\tau'$, so that $\tau':=exp(ad_\theta)(d+\tau)-d$. We want to show that $\tau'\in \Theta_{-1}$. Clearly $\tau'$ reduces topological degree by $1$, as $d+\tau$ did. So we now need to show that $\tau'$ reduces resolution degree by $\geq 2$. But 
$$\tau'=[\sum_{i=0}^\infty ad_\theta^i(d+\tau)]-d=[\sum_{i=1}^\infty ad_\theta^i(d+\tau)]+\tau$$
As $\theta\in\Theta_0$, $ad_\theta^i$ reduces resolution degree by at least $1$ for all $i>0$. So do $d$ and $\tau$, so $ad_\theta^i(d+\tau)$ reduces resolution degree by at least $2$ for $i>0$. $\tau$ also reduces resolution degree by at least $2$ as $\tau\in\Theta_{-1}$. So $\tau'$ reduces resolution degree by at least $2$. Thus $\tau'\in\Theta_{-1}$.

\item[b)] Now we wish to show that $(d+\tau)^2=0 \implies (d+\tau')^2=0$ for $\theta$, $\tau$, and $\tau'$ as above. Recall that $$(d+\tau')^2=[d,\tau']+\frac{1}{2}[\tau', \tau']$$
 We first consider $[\tau',\tau']$. Note that $$[\tau',\tau']=[exp(ad_\theta)(d+\tau)-d,exp(ad_\theta)(d+\tau)-d]$$ $$=[exp(ad_\theta)(d+\tau),exp(ad_\theta)(d+\tau)]-2[exp(ad_\theta)(d+\tau),d]+[d,d]$$
 We know that $[d,d]=0$. Furthermore, as $\theta$ is a locally nilpotent derivation, we have that $exp$ preserves the bracket, and so $[exp(ad_\theta)(d+\tau),exp(ad_\theta)(d+\tau)]=exp(ad_\theta)[d+\tau,d+\tau]$. But $(d+\tau)^2=0$ by assumption, so $[exp(ad_\theta)(d+\tau),exp(ad_\theta)(d+\tau)]=0$. Thus our equation reduces to $$[d,exp(ad_\theta)(d+\tau)-d]+\frac{1}{2}(-2[exp(ad_\theta)(d+\tau),d])$$
 
However, $$[d,exp(ad_\theta)(d+\tau)-d]=[d,exp(ad_\theta)(d+\tau)]$$$$=(-1)^{(-1)(-1)}[exp(ad_\theta)(d+\tau),d]=[exp(ad_\theta)(d+\tau),d]$$

So, $$[d,exp(ad_\theta)(d+\tau)-d]+\frac{1}{2}(-2[exp(ad_\theta)(d+\tau),d])$$$$=[exp(ad_\theta)(d+\tau),d]-[exp(ad_\theta)(d+\tau),d]=0$$

Thus $(d+\tau')^2=0$, as desired.

It follows that $\tau'$ is a perturbation. $exp(ad_\theta )$ is, then, the weak equivalence required.
\end{description}

\underline{ $\implies$:}

Assume there exists a quasi-isomorphism $\phi :(\mathcal{L},d+\tau,i_{\tau})\to (\mathcal{L}, d+\tau',i_{\tau'} )$. We want to show that $(\tau',i_{\tau}') = exp(ad_\theta )(\tau,i_{\tau})$ for some $\theta \in \Theta_0 $.

\

We will construct a unipotent, homology-preserving DGLA homomorphism $\Phi:(\mathcal{L},d+\tau)\to (\mathcal{L}, d+\tau')$ by induction on resolution degree in the bigraded model. 

\begin{enumerate}
\item For $x\in \mathcal{L}_*^0$, define $\Phi x=x$. Note that as this is the identity, it is clearly a DGLA homomorphism and unipotent. As all homology classes have a representative in resolution degree $0$, this will ensure that $\Phi$ preserves homology classes at all levels as long as it remains to be a DGLA homomorphism.

\

\item For $y\in \mathcal{L}_*^1$, define $\Phi y =y$ again. As above, $\Phi$ is still a unipotent DGLA homomorphism, as $d+\tau=d=d+\tau'$ in this degree.

\

\item Let $z\in \mathcal{L}_*^2$. We know $\Phi (d+\tau) z=(d+\tau)z$, and so as $\Phi$ is a DGLA homomorphism, $(d+\tau)z$ is a boundary in $(\mathcal{L},d+\tau')$. So, $(d+\tau)z-(d+\tau')z=(\tau-\tau' )z$ is as well, and thus $(\tau-\tau' )z=(d+\tau')w$ for some $w$. It remains to show that there is such a $w$ with $mrd(w)<2$. Note that $[dz]_{d+\tau}=[-\tau]_{d+\tau}$, so $[\Phi dz]_{d+\tau'}=[\Phi (-\tau)]_{d+\tau'}$. But as $\Phi$ is the identity in these resolution degrees, we get that $[dz]_{d+\tau'}=[-\tau]_{d+\tau'}$. Furthermore, $[dz]_{d+\tau'}=[-\tau']_{d+\tau'}$. This gives us that $[-\tau]_{d+\tau'}=[-\tau']_{d+\tau'}$. As $\tau$ and $\tau'$ are both in resolution degree $0$ and represent the same class in homology, by the construction of the minimal model we have that there exists an element $w$ such that $dw=\tau-\tau'$. This $w$ is in resolution degree $1$, as desired.

So, for all $z\in \mathcal{L}_*^2$ we define $\Phi z=z+w$, where $w$ is found for each $z$ as above. Clearly $\Phi$ is unipotent by construction, so we need only confirm that $\Phi$ remains a DGLA homomorphism by the end of the construction process. But $(d+\tau')\Phi z=(d+\tau')(z+w)=(d+\tau')z+(\tau-\tau')z=(d+\tau)z=\Phi(d+\tau)z$, and so $\Phi$ is still a DGLA homomorphism.

\

\item We now proceed to the inductive step. Let $s\in \mathcal{L}_*^n$ and assume that $\forall z\in \mathcal{L}_*^{<n}$ we have that $\Phi z=z+w$, where $mrd(w)<mrd(z)$, and that $\Phi$ satisfies both the unipotence and the DGLA homomorphism conditions when restricted to resolution degrees $<n$. Then $\Phi (d+\tau)s=(d+\tau)s+w$ is a boundary in $(\mathcal{L},d+\tau')$ and so there exists a $v\in \mathcal{L}$ such that $(d+\tau')v=(d+\tau)s+w-(d+\tau')s=(\tau-\tau')s+w$. Write $v=\Sigma_{k=0}^m v_k$, where $v_k$ is a homogeneous element of resolution degree $k$. If $m<n$ we're done, so assume $m\geq n$. Then $d v_m=0$ as $mrd((\tau-\tau')s+w)\leq n-2$ and $resdeg(d v_m)=m-1\geq n-1$. As $d v_m=0$ and $v_m$ is in resolution degree $m>0$, $\exists u_m$ s.t $du_m=v_m$. We define $\tilde{v}:=(\Sigma_{k=0}^{m-1} v_k)-\tau u_m=\Sigma_{k=0}^{m-1} \tilde{v}_k$, which has the property $mrd(\tilde{v})<mrd(v)$. If $m-1<n $ we're done. If not, continue replacing the top resolution degree piece of $\tilde{v}$ in the same manner and redefining $\tilde{v}$ until $mrd(\tilde{v})\leq n-1$. Now we at least have an element of the correct resolution degree.

However, we have to ensure that our replacing process does not change the element's image under $d+\tau'$. That is, we want $(d+\tau')v=(d+\tau')\tilde{v}$. It suffices to show that $(d+\tau')(-\tau u_m)=(d+\tau')v_m$ for any pair of elements $u_m$, $v_m$ with $d u_m=v_m$. But $0=(d+\tau')^2u_m=(d+\tau')(d+\tau')u_m=(d+\tau')(v_m+\tau' u_m)\implies (d+\tau')(-\tau u_m)=(d+\tau')v_m$, as desired. So, $mrd(\tilde{v})<n$ and $(d+\tau')\tilde{v}=(\tau - \tau')s+w$. Define $\Phi s=s+\tilde{v}$. We must show that $\Phi$ remains to be a DGLA homomorphism. But $(d+\tau')\Phi s=(d+\tau')(s+\tilde{v})=(d+\tau')s+(\tau - \tau')s+w=(d+\tau )s+w=\Phi (d+\tau )s$, as desired.

\

\item Now that we have produced a unipotent DGLA homomorphism $\Phi:(\mathcal{L},d+\tau)\to (\mathcal{L}, d+\tau')$, we still need to find a $\theta \in \Theta_0$ such that $d+\tau'=exp(\theta)(d+\tau)$. However, $\theta:=log(\Phi - id)$ satisfies the conditions to be in $\Theta_0$ and, since $exp(ad_{Log(\Phi - id)})=\Phi$ and $\Phi:d+\tau\to d+\tau'$, we have actually found the $\theta$ necessary. 
\qed
\end{enumerate}

\begin{corollary}In the sense of triples, two perturbations $\tau$ and $\tau'$ of a bigraded model $(\mathcal{L},d)$ are equivalent if and only if they differ by an action of $exp$. That is, $(d+\tau')=\displaystyle \prod_{k=1}^n(exp(ad_{\theta_k})(d+\tau))$ if and only if $(\mathcal{L},d+\tau,i_{\tau})$ and $(\mathcal{L},d+\tau',i_{\tau'})$ are quasi-isomorphism equivalent.
\end{corollary}

\emph{Proof:}

\underline{$\implies$}: If $d+\tau'$ is a product of $exp(ad_{\theta_k})(d+\tau,i_{\tau})$s, then each one gives us an arrow, and the chain of arrows gives us the equivalence.

\underline{$\impliedby$}: If $(\mathcal{L},d+\tau,i_{\tau})\sim (\mathcal{L}, d+\tau',i_{\tau'} )$ then we have a chain of weak equivalences from $(\mathcal{L},d+\tau,i_{\tau})$ to $(\mathcal{L},d+\tau',i_{\tau'})$, potentially with both left- and right-facing arrows. Each map gives rise to a $\theta$ from the theorem, and these $\theta$s give us $d+\tau'$ under the image of $exp$, with some $exp(\theta)$s needing to be inverted when the corresponding equivalence arrow was in the wrong direction. 
\qed

\

With the results above, the moduli space $\mathcal{M}_P$ for all rational homotopy types of simply connected topological spaces with homotopy Lie algebra $P$ can be constructed as follows:

\begin{enumerate}
\item Begin with a homotopy Lie algebra $P$ and the corresponding bigraded model $(\mathcal{L}_P,d)$.
\item Consider the variety $V$ of deformations $\tau \in \Theta_{-1}$ such that $(d+\tau)^2=0$.
\item Mod out by $exp(ad_{\Theta_{0}})$.
\item Mod out by $Aut(P)$.
\end{enumerate}

\subsection{A Technical Summary}

We first begin with a homotopy Lie algebra $P$. As $P$ is a graded Lie algebra, we can choose a graded vector space $V$ of indecomposibles that generate $P$. From this vector space we can begin to build the bigraded model. Specifically, we begin with the first level, resolution degree $0$, of the bigraded model by allowing $V$ to generate a free graded Lie algebra $FV$ with differential $0$. However, the homology of $(FV,0)$ is $FV$, which is not necessarily $P$. If $FV\neq P$, then this is not a bigraded model for $P$. To rectify this, we choose representatives for the extra homology classes in $FV$ and introduce generators in resolution degree $1$ together with differential maps to kill those extra homology classes. Then we generate a free bigraded DGLA. The homology of this bigraded DGLA is then checked against $P$. The process continues until there no longer is any inconsistencies between the homology of the DGLA and the original homotopy Lie algebra $P$. At the end, we get the bigraded model for $P$.

This bigraded model comes with a map $i$ from its homology to $P$. As a consequence of Lemma 3.19, this map is actually induced from the map $\mathcal{L}^0 \to P$, and thus a map $\hat{i}:(\mathcal{L},d)\to P$.
$$
\bfig
\morphism(0,0)|a|/{->}/<0,-500>[\mathcal{L}_0`P;\hat{i} ]
\morphism(750,0)|a|/{->}/<0,-500>[H(\mathcal{L},d)`P;i]

\efig
$$

However, there isn't necessarily a map from the DGLA to $P$ after we perturb the bigraded model. In particular, $(\mathcal{L},d+\tau)$ would only have a quasi-isomorphism to $P$ if the perturbation itself was trivial. So for a nontrivial perturbation we don't have such a map. However, there is still necessarily a map on homology.

$$
\bfig
\morphism(0,0)|a|/{-->}/<0,-500>[(\mathcal{L},d+\tau)`P;? ]
\morphism(750,0)|a|/{->}/<0,-500>[H(\mathcal{L},d+\tau)`P;i_\tau ]

\efig
$$
This shifts our focus to, in fact, triples $(\mathcal{L},d+\tau, i_\tau)$ as the maps to $P$ ensure that the models are of the correct rational homotopy type.

We have the bigraded model and perturbations established, leading to a collection $\{(\mathcal{L},d+\tau,i_\tau)|i_\tau : H(\mathcal{L},d+\tau)\to P \textrm{\ is a GLA isomorphism}\}$. We want to know when two perturbations actually represent the same rational homotopy type. For two perturbations to represent the same rational homotopy type, they would need to differ by a finite zig-zag of quasi-isomorphisms $$(\mathcal{L},d+\tau,i)\to (\mathcal{L},d+\tau_1,i_1)\leftarrow (\mathcal{L},d+\tau_2,i_2) \to ... \leftarrow (\mathcal{L},d+\tau_k,i_k) \to (\mathcal{L},d+\tau',i_{\tau '})$$ where each quasi-isomorphism is an action of $exp(ad_{\Theta_0})$, potentially adjusted by a lift of an automorphism of $P$ bringing $i_j\to i_{j+1}$. Note in particular that each $exp(ad_{\theta})$ respects the homology isomorphisms $i_j\to P$. Then, we need only mod out by automorphisms of $P$ to ensure that each rational homotopy type is only counted once. So, we produce $\mathcal{M}_P$ as: $$\big\{\{(\mathcal{L},d+\tau,i_\tau)|i_\tau : H(\mathcal{L},d+\tau)\to P \textrm{\ is a GLA isomorphism}\}/exp(ad_{\Theta_0})\big\}/Aut(P)$$
Alternatively, one could view the construction of $\mathcal{M}_P$ by considering all possible $(\mathcal{L},d+\tau,i_\tau)$, with $i_\tau$ a GLA isomorphism as above, and then modding out by the group action induced by $exp(ad_{\Theta_0})\times Aut(P)$. Either way, we start with all possible perturbations that would yield a differential (in terms of topological degree) on the bigraded model. Then we mod out by things that would give rise to two equivalent perturbations somehow. The first idea does this stepwise, while the second reduces to a single group action by the product group $exp(ad_{\Theta_0})\times Aut(P)$.

\subsection{Examples of Moduli Spaces}

We close with some brief examples of the moduli spaces described above, to give the reader some insight into the various amounts of complexity such spaces might have.

\begin{example}[$S^2$]
Recall from Example 3.10 the bigraded model of $S^2$:

\begin{center}

\begin{tabular}{c c c c c c}
0 & $a$ & $[a,a]$ & $0$ & $0$ & ... \\
\hline
\hline
0 & 1 & 2 & 3 & 4 & ... \\
\end{tabular}
\end{center}
Since there are no elements in positive resolution degree, there is no room for perturbations to occur. Thus the moduli space for rational homotopy types with the same homotopy Lie algebra is just a single point. Because of this rather special property, $S^2$ is called \textit{intrinsically coformal}.
\end{example}

\begin{example}[$\mathbb{C}P^2$]
Recall that the bigraded model for the coformal space with the same homotopy Lie algebra as $\mathbb{C}P^2$ is 

\begin{center}
\begin{tabular}{c c c c c c c}
0 &  0  &   0     &  0  &    0     &       $c$     & $[b,b];[c,a]$ \\
0 &  0  &   0     & $b$ & $[b,a]$  & $[b,[a,a]]$ & $y$  \\
0 & $a$ & $[a,a]$ & 0 & $x$ & $[x,a]$ & ... \\
\hline
\hline
0 & 1 & 2 & 3 & 4 & 5 & ... \\
\end{tabular}
\end{center}
with $db=[a,a]$ and $dc=[b,a]$. There is a nontrivial perturbation $\tau$ with $\tau c=x$, which represents the rational homotopy type of $\mathbb{C}P^2$.  In \cite{NM}, they show that this moduli space is in fact a two-point space by exporting the problem to the DGA category. 

\end{example}

\begin{remark}
The example above begins to touch upon the idea of \textit{Massey brackets}, an idea tightly related to the variability of these moduli spaces. In particular, $[b,a]$ is the triple bracket $[a,a,a]$. A full treatment of Massey brackets is beyond the scope of this paper, but the interested reader can refer to \cite{SS3}.
\end{remark}

\begin{example}[An Interesting Moduli Space]
Let $P$ be the following graded Lie algebra:

\begin{center}
\begin{tabular}{c c c c c c c}
0 & $a, b$ & $[a,a]$ & $c$ & $[c,a];[c,b];e$ & ... & ... \\
\hline
\hline
0 & 1 & 2 & 3 & 4 & 5 & ... \\
\end{tabular}
\end{center}
Then $(\mathcal{L}_P,d)$ is the following (in low topological degree, with some elements omitted due to spatial constraints):

\begin{center}
\begin{tabular}{c c c c c c c}
 &       &       &        &                          & $w;z$ \\
 &       &        & $x;y$ & $[a,x];[b,x];[a,y];[b,y]$& ... \\ 
0 & $a, b$ & $[a,a];[b,b];[a,b]$ & $c;[a,[b,b]];[a,[a,b]]$ & $[c,a];[c,b];e;[[a,a],[b,b]];...$ & ... \\
\hline
\hline
0 & 1 & 2 & 3 & 4 & 5 \\
\end{tabular}
\end{center}
with differential as defined below:
$$dx=[b,b]$$
$$dy=[a,b]$$
$$dw=[b,x]$$
$$dz=[a,x]+2[b,y]$$
Already we can see the model becoming incredibly complex as topological and resolution degrees increase. The possible nontrivial perturbations of the differential on $w$ and $z$ are as follows:
$$\tau w=e$$
$$\tau w=[c,a]$$
$$\tau w=[c,b]$$
$$\tau z=e$$
$$\tau z=[c,a]$$
$$\tau z=[c,b]$$
None of these are equivalent to each other, and linear combinations of these can also be perturbations of the differential. There are also some trivial perturbations, as $[[a,a],[b,b]]$, $[[a,a],[a,b]]$, and $[[a,b],[b,b]]$ would all be acceptable targets for perturbations but are in the image of $d$ (the corresponding elements in bidegree ($5,1$) are not difficult to determine). So, even only up to topological degree $5$, this moduli space has a multitude of parameters to exploit.
\end{example}

\end{document}